\documentstyle[11pt]{article}
\begin{document}
\bibliographystyle{plain}
\newtheorem{prob}{Problem}[section]
\newtheorem{theo}{Theorem}[section]
\newtheorem{proper}{Property}[section]
\newtheorem{lem}{Lemma}[section]
\newtheorem{deff}{Definition}[section]
\newtheorem{cor}{Corollary}[section]
\newtheorem{pro}{Proposition}[section]
\newtheorem{conj}{Conjecture}[section]
\newenvironment{proof}{\noindent {\bf
Proof.}}{\rule{3mm}{3mm}\par\medskip}
\title{OTIS Layouts of De Bruijn Digraphs
\thanks{Supported by the NNSFC Grant No. 19971056
and the RFDP Grant No. 2000024837.} }
\author{Yaokun Wu, Aiping Deng\\
{\small Department of Mathematics, Shanghai Jiao Tong University,}\\
{\small 1954 Huashan Road, Shanghai, 200030, China}\\
{\small Email: ykwu@sjtu.edu.cn,
apdeng@sjtu.edu.cn}\\
 \date{}}
\maketitle

\noindent {{{\bf Abstract}}

The Optical Transpose Interconnection System (OTIS) was proposed
by Marsden et al. [Opt. Lett 18 (1993) 1083--1085]
 to implement very dense one-to-one interconnection between
 processors in a free space of optical interconnections.
The system which allows one-to-one optical communications
 from $p$ groups of $q$ transmitters to $q$ groups of $p$ receivers, using
 electronic  intragroup communications for each group of
 consecutive $d$ processors,
  is denoted by  OTIS$(p,q,d)$.
 $H(p,q,d)$ is the digraph which characterizes the underlying
 topology of the optical interconnection implemented by OTIS$(p,q,d)$.
 A digraph has an OTIS$(p,q,d)$ layout if it is isomorphic to
  $H(p,q,d).$ Based on results of
 Coudert et al. [Networks 40 (2002) 155--164], we show that  De Bruijn digraph
 $B(d,n)$ has an OTIS$(d^{p'},d^{n+1-p'},d)$ layout if and only if $\gcd (p',n+1)=1.$
 We also prove that  $H(p,q,d)$ is a line digraph if and only if $d\mid \gcd (p,q).$

 }
\section{INTRODUCTION}
\noindent It is known that electronic interconnects perform better
than optical interconnects when the distance is up
 to
 a few millimeters \cite{feld88}, while on a longer distance,
the latter has many  advantages over the former, like less
crosstalk,
 less power consumption, higher speed, and
 high bandwidth channels at a single communication point \cite{day02,zane00}.
 To take benefits from both optical and  electronic technologies, Marsden et al. \cite{mar93} proposed the
Optical Transpose Interconnection System (OTIS), which has gained
considerable attention
\cite{chi96,C&F&P02,C&F&M1,coudert00,C&F&M00,day02,wang98,zane00}.
The idea of the OTIS system is to partition the processors  into
groups and use electronic interconnects for the intragroup
communications(which are of short distance) while optical
interconnects for the intergroup communications  (which are of
larger distance). For the purpose of designing a well-behaved OTIS
system, we may hope to use some good topological structures as
models for the electronic interconnects as well as the optical
interconnects. Since arbitrary connections using optical links via
lenses are harder to implement than using wires on a VLSI
 circuit, multi-chip module or printed circuit board, how to realize a given
 good topology as optical
interconnects has been of special concern.

After the brief introduction of the background, we now turn to a
mathematical abstraction of the OTIS layout problem. Let us follow
the model of Coudert et al. \cite{C&F&P02}. For any two integers
$a\leq b$, write $[a,b]$ for the set $\{ c\in Z \mid a\leq c \leq
b \}.$ Suppose $p,q$ are two positive integers. Any $a\in
[0,pq-1]$ can be uniquely expressed as $a=iq+j$, where $j\in
[0,q-1]$ and $i\in[0,p-1].$ We use the notation
 $(i,j)_{p,q}$ for such an $a$ and say that $a$ has first
$(p,q)$-coordinate $i$ and second $(p,q)$-coordinate $j$.
 Let there
be $pq$ processors,  each equipped with an optical
transmitter/receiver pair. We label these processors with
$[0,pq-1]$. Set two planes of lenses arrays, one consisting of $p$
lenses each corresponding to $q$
 transmitters for a group of processors with the same first $(p,q)$-coordinate,
  the other consisting of $q$ lenses each corresponding to  $p$
 receivers for a group of processors with the same first $(q,p)$-coordinate.
 These lenses  establish optical links from transmitters
 of processor $(i,j)_{p,q}$ to receivers of processor $(q-1-j,p-1-i)_{q,p}$
 for $i\in[0,p-1],j\in[0,q-1]$,
 namely the directed optical links
 transpose both coordinates. For a divisor $d>1$ of $pq,$ we
 use electronic interconnects among processors with labels $[kd,kd+d-1]$
for $k\in[0,pq/d-1].$ This way, we then have built an
OTIS$(p,q,d)$ architecture. Viewing each group of processors
$[kd,kd+d-1]$ as a node and assigning as many  arcs from node
$[kd,kd+d-1]$ to $[k'd,k'd+d-1]$ as there are  optical links from
transmitters of the first group to receivers of the second, we
obtain a digraph $H(p,q,d)$, which reflects the optical intergroup
communication pattern of the OTIS$(p,q,d)$ architecture. Note that
we can also use two sets of processors in the above construction,
one corresponding to the transmitters and the other the receivers.
The basic connecting unit thus obtained may be cascaded  to
accommodate successive processing planes. It is not hard to see
that the optical intergroup communication pattern of the resulting
structure is a multistage interconnection network \cite{Feng} and
$H(p,q,d)$  characterizes its underlying topology. Observe that
$H(p,q,d)$ is a $d$-regular digraph on $pq/d$ vertices. We now
come to
\begin{deff}  (\cite{C&F&P02} Definition 4.2)
We say that a $d$-regular digraph $G$ has an OTIS layout provided
there are positive integers $p,q$ such that $G$ is isomorphic to
$H(p,q,d)$. \end{deff}
The OTIS layout problem is the problem to
characterize all OTIS layouts for a given digraph and to find
among all OTIS layouts the one which is optimal in some aspects,
like using the fewest lenses, namely minimizing $p+q.$

The technique of line digraph iterations proves to  be useful  in
producing vast families of good network models
\cite{b0,cao99,fiol84,wu03}. Particularly, for any positive
integers $d$ and $n,$ the $n$th line digraph of the complete
digraph on $d$ vertices with loops, $L^n(K_d^+)$, called the
$n-$dimensional $d-$ary {\it De Bruijn digraph} and denoted by
$B(d,n)$, has been the focus of much study as a very good
interconnection
structure\cite{B2,B3,C&F&P02,C&F&M1,C&F&M00,H&H01,litman,S&H&F00}.
It  thus seems natural to  address the isomorphisms between
$H(p,q,d)$ and De Bruijn digraphs, or more generally, line
digraphs.

Coudert et al. deduced the following characterization of
OTIS$(d^{p'},d^{q'})$ layouts for De Bruijn digraphs.

\begin{theo} (\cite{C&F&P02} Lemma 4.4) Let $p'+q'-1=n.$ For any degree
$d$,  $B(d,n)$ and $H(d^{p'},d^{q'},d)$ are isomorphic if and only
if the permutation $f$ of $Z_{n}$ defined by
$$f(i)=\left\{\begin{array}{ll}i+p' &  if\,\,
i\in\{0,1,\cdots,q'-2\};\\
p'-1 & if\,\,
i=q'-1;\\i+p'-1 \pmod{n} & otherwise,\\
\end{array}
\right.$$is cyclic. \label{codert}
\end{theo}

As a corollary, Coudert et al. (\cite{C&F&P02} Corollary 4.8)
pointed out that whether or not $B(d,n)$ has an
OTIS$(d^{p'},d^{q'})$ layout can be checked in $O(n)$ time.
Finally, they concluded \cite{C&F&P02} by indicating that
 their exhaustive search led them to

\begin{conj}
If $B(d,n)$ has an OTIS$(p,q)$ layout, then $p,q$ must be powers
of $d.$ \label{conj}
\end{conj}

Our paper is an effort to characterize all OTIS layouts of De
Bruijn digraphs. Making use of Theorem \ref{codert}, we  will show
in Section \ref{2} that  De Bruijn digraph
 $B(d,n)$ has an OTIS$(d^{p'},d^{n+1-p'})$ layout if and only if $\gcd (p',n+1)=1.$
 Note that using Euclidean algorithm, we  only need $O(\log n)$ time steps
 to evaluate $\gcd (p',n-1)$ (\cite{wilf94} Theorem 4.2.1) and thus it implies an improvement of the above-mentioned
 result of Coudert et al. from $O(n)$ to  $O(\log n).$ As a step toward proving the conjecture of
 Coudert et al., we will prove in Section \ref{3} that $H(p,q,d)$
 is a line digraph if and only if $d\mid \gcd (p,q).$

\section{DE BRUIJN DIGRAPH}
\label{2} \noindent
 Let $p',q',n$
 be three positive integers such that $p'+q'-1=n.$
 Define a permutation $g_{p',q'}$ on $[0,n-1]$
by
$$g_{p',q'}(i)=\left\{\begin{array}{ll}i+p' &  \mbox{if}\,\,
i\in [0,q'-2];\\
i+p'-q' & \mbox{if}\,\,
i=q'-1;\\i-q' & \mbox{if}\,\,i\in [q',n-1].\\
\end{array}
\right.$$ We will adopt the convenient notation $g$ for
$g_{p',q'}$ hereafter. $g$ is just another representation of the
$f$ as defined in Theorem \ref{codert}.  Thus we have

\begin{theo}$B(d,n)$ and $H(d^{p'},d^{q'},d)$ are isomorphic
if and only if $g$ is a cyclic permutation on  $[0,n-1].$
\label{re}
\end{theo}

As we will see immediately, the form of $g$ is more suitable for
an investigation of its cycle structure and the above trivial
reformulation of Theorem \ref{codert} is indeed a key observation
for us.
 Let $\lambda =\gcd (p',q').$ For each $i\in [0,n-1],$
write $C_i$ for the set $\{j\in [0,n-1]\mid j\equiv i
\pmod{\lambda} \}$ and $O_i$ for the orbit of $i$ under the action
of $g.$ We use $\mid S\mid$ for the cardinality of any finite set
$S.$ The Kronecker Delta $\delta _{i,j}$ is defined as having
value $1$ when $i=j$ and $0$ otherwise.
\begin{theo} $g$ has exactly $\lambda$ orbits. Indeed, the partition
of $[0,n-1]$ into orbits of $g$ is the same with its partition
into congruence classes modulo $\lambda ,$ namely $O_i=C_i$, $i\in
[0,n-1].$ \label{orbit}
\end{theo}
\begin{proof}
Clearly, it always holds $g(i)-i \equiv 0 \pmod {\lambda}.$ This
means $O_i\subseteq C_i$ for all $i\in [0,n-1].$ But $[0,n-1]$ is
a disjoint union of $C_i$ for $i\in S=[0,\lambda -2]\cup
\{q'-1\}.$ Moreover, it is easy to see that $\mid C_i\mid
=\frac{n+1}{\lambda}-\delta _{i,q'-1}$ for $i\in S.$ Thus our goal
is just to verify that $\mid O_i\mid \geq
\frac{n+1}{\lambda}-\delta _{i,q'-1}$ for $i\in S.$

Take any $i\in S.$ Let $\alpha _i=\mid O_i\cap [0,q'-2]\mid,$
$\beta _i=\mid O_i\cap \{q'-1\}\mid,$ and $\gamma _i=\mid O_i\cap
[q',n-1]\mid .$  As $O_i$ is an orbit, we obtain $0=\sum _{j\in
O_i}j-\sum _{j\in O_i}j=\sum _{j\in O_i}(g(j)-j)=\alpha _ip'+\beta
_i (p'-q')-\gamma _i q'=(\alpha _i+\beta _i)p'-(\beta _i+\gamma
_i)q'. $ Cancelling the common factor $\lambda$ of $p'$ and $q'$
yields  $(\alpha _i+\beta _i)\frac{p'}{\lambda}=(\beta _i+\gamma
_i)\frac{q'}{\lambda}.$
 Since $\gcd (p'/\lambda ,q'/\lambda)=1,$ it follows that
$\alpha _i+\beta _i$ is a multiple of $q'/\lambda$ and therefore
$q'/\lambda \leq
 \alpha _i+\beta _i.$ Similarly, we have $p'/\lambda \leq
 \beta _i+\gamma _i.$ These two inequalities together implies
 that $\alpha _i+\beta _i+\gamma _i\geq \frac{p'+q'}{\lambda}-\beta _i=\frac{n+1}{\lambda}-\beta _i.$
But it holds $\alpha _i+\beta _i+\gamma _i=\mid O_i\mid.$
Furthermore, we can derive from $O_i \subseteq C_i$  that $\beta
_i=\delta _{i,q'-1}.$ So we have arrived at $\mid O_i\mid \geq
\frac{n+1}{\lambda}-\delta _{i,q'-1},$ as desired.
\end{proof}

Notice that $\gcd (p',q')=\gcd (p',p'+q')=\gcd (p',n+1)$.
Consequently, by Theorem \ref{re} and Theorem \ref{orbit}, we can
establish the following characterization of OTIS layouts of De
Bruijn digraphs.
\begin{theo} For $p'\in [0,n+1],$ $B(d,n)$ and $H(d^{p'},d^{n+1-p'},d)$ are isomorphic
if and only if $\gcd (p',n+1)=1.$ \label{final}
\end{theo}

\section{LINE DIGRAPH}
\label{3} \noindent We remark that, assuming Conjecture
\ref{conj}, which holds trivially  when $d$ is a prime, Theorem
\ref{final} tells us that there are totally $\phi (n+1)$ different
OTIS layouts for $B(d,n),$ where $\phi$ is the Euler's totient
function.
 But is Conjecture \ref{conj} really true in general cases?
 As a prominent characteristic of De Bruijn digraphs is their iterated line digraph
 structure
 \cite{wu03}, we are naturally led to the study of those parameters
 $p$ and $q,$ such that for a fixed $n,$
 $H(p,q,d)$ is an $n$th iterated line digraph.
This line of research requires some preliminary results on
characterizing iterated line digraphs. A classic result is
Heuchenne's characterization of line digraphs
\cite{berge,Hem,Heuch}, proved about 40 years ago. Indeed, our
subsequent work on characterizing OTIS layouts of line digraphs is
just based on it.  For possible later use in tackling the problem
for general iterated line digraphs, instead of merely presenting
Heuchenne's characterization, we include here a characterization
of iterated line digraphs, which generalizes Heuchenne's result
and an earlier generalization of it due to Beineke and Zamfirescu
\cite{Bei}.

 A digraph is said to
satisfy the $n$th {\em Heuchenne condition} \cite{Bei,Heuch} if
for any of its vertices $u$, $v$, $w$, and $x$
 (not necessarily distinct)
for which there exist $n$-walks from $u$ to $w$, from $v$ to $w$,
and from $v$ to $x$, there must also exist an $n$-walk from $u$ to
$x.$  Restricting our attention to the case of $n=1,$ the
following theorem  is just  Heuchenne's characterization.

\begin{theo} (\cite{wu03} Theorem 7) \label{line}
Let $G$ be a digraph without sinks or sources. Then $G$ is an
$n$th  line digraph if and only if the following conditions are
satisfied:
\begin{enumerate}
\item[{\rm ($I$)}]
There are no multiple $n$-walks between any pair of vertices;
\item[{\rm ($II$)}]
$G$ satisfies both the  $(n-1)$th  and the $n$th Heuchenne
conditions.
\end{enumerate}
\end{theo}

 We are in a position
to prove our main results for the OTIS layouts of line digraphs.
For $j\in [id,id+d-1],$ we define $\alpha (j)=i,$ namely $\alpha
(j)=\lfloor\frac{j}{d}\rfloor.$

\begin{theo}\label{lem}
If $d\mid \gcd (p,q),$ then $H(p,q,d)$ is a line digraph.
\end{theo}
\begin{proof} For any $i\in [pq/d],$ it can be uniquely expressed as
$i=(t,s)_{p,\frac{q}{d}}.$ Let $v_i$ be the vertex of $H(p,q,d)$
corresponding to the interval $M_i=[di,di+d-1].$ It is
straightforward to check that $M_i=\{(t,ds)_{p,q},\cdots
,(t,ds+d-1)_{p,q}\}.$ Thus the out-neighbors of $v_i$ can be
enumerated as $v_{i_0},\cdots ,v_{i_{d-1}}$ such that
$(q-1-ds,p-1-t)_{q,p}\in M_{i_0},\cdots
,(q-1-(ds+d-1),p-1-t)_{q,p}\in M_{i_{d-1}}.$ It follows that the
out-neighbors of $v_i$ are just $v_{\alpha
((q-1-ds,p-1-t)_{q,p})},\cdots, v_{\alpha
((q-1-(ds+d-1),p-1-t)_{q,p})},$ which turns out to be
 $v_{(q-1-ds,\alpha (p-1-t))_{q,\frac{p}{d}}},\cdots, v_{(q-1-(ds+d-1),\alpha (p-1-t))_{q,\frac{p}{d}}}.$
Since these $d$ vertices are obviously pairwise different,
$H(p,q,d)$ fulfils condition (I). Also, we see that for any two
vertices $v_i$ and $v_{i'}$ with $i=(t,s)_{p,\frac{q}{d}},$ and
$i'=(t',s')_{p,\frac{q}{d}},$ respectively, their out-neighbor set
will be disjoint if $s\not= s'$ or $\alpha (p-1-t)\not= \alpha
(p-1-t'),$ and will be identical otherwise. This shows that
condition (II) holds as well. By Heuchenne's characterization,
this then completes the proof.\end{proof}

Interestingly, the converse of Theorem \ref{lem} is also true,
which  provides  partial support to Conjecture \ref{conj}. Recall
that for a digraph $G$, its dual, written $\overleftarrow{G},$ is
the digraph obtained from $G$ by reorienting each edge in the
opposite direction as in $G.$
\begin{theo}\label{theo}
If $H(p,q,d)$ is a line digraph, then $\gcd (p,q)$ is a multiple
of $d.$
\end{theo}

\begin{proof} Clearly, it holds $H(p,q,d)=\overleftarrow{H(q,p,d)}$
and $L(\overleftarrow{G})=\overleftarrow{L(G)}$ for any digraph
$G.$ Thus we only need to prove $d \mid p.$ Let us fix some
notation before proceeding. Hereafter, the digraph $H(p,q,d)$ is
simply called $H.$ For any $i,$ we denote by $I_i$ the interval $[
\alpha ((i,0)_{p,q})d, \alpha ((i,0)_{p,q})d+d-1]$  and by $v_i$
the vertex of $H$ which corresponds to $[id,id+d-1].$

We first claim that $d\leq \min (p,q).$ Again, we only need to
prove $d\leq p$ due to the fact
$L(\overleftarrow{G})=\overleftarrow{L(G)}.$ Seeking a
contradiction, suppose that $d>p.$ Then, both $(0,0)_{q,p}$ and
$(1,0)_{q,p}$ belong to $[0,d-1],$ the interval corresponding to
$v_0.$ But $[pq-d,pq-1]$ includes $(p-1,q-1)_{p,q}$ and
$(p-1,q-2)_{p,q},$ which send links to $(0,0)_{q,p}$ and
$(1,0)_{q,p},$ respectively. Hence there are multiple arcs from
$v_{\frac{pq}{d}-1}$ to $v_0,$ in violation of condition (I).

Next  observe that condition (II) means that for the line digraph
$H$ the relation $\sim$ of having a common out-neighbor  is an
equivalence relation on its vertex set. Because $H$ is
$d$-regular, each equivalence class of $\sim$ has size $d.$
 Moreover, since
the vertices  of $H$  correspond to pairwise disjoint intervals of
length $d,$ the relation $\sim$ naturally extends to the
equivalence relation $\approx$ on $[0,pq-1]$ such that $i\approx
j$ if and only if $v_{\alpha (i)}\sim v_{\alpha (j)}.$ Each
equivalence class under $\approx$ has a size $d$ times as large as
that of an equivalence class of $\sim$ and thus contains $d^2$
elements.

To finish the proof, suppose, on the contrary,  that $p$ is not a
multiple of $d.$ Note that $p\geq d$ and $d\mid pq.$ It
immediately follows that $p>d$ and $\alpha ((q-1,0)_{q,p})=\alpha
((q-2,p-1)_{q,p}).$
  Because there are links from $(p-1,0)_{p,q}$ to $(q-1,0)_{q,p}$
 and from $(0,1)_{p,q}$ to $(q-2,p-1)_{q,p},$
 the latter formula tells us that
 \begin{equation}(p-1,0)_{p,q}\approx (0,1)_{p,q}.
 \label{final}
 \end{equation}
Further notice that for each $i\in [0,d-1],$ $(i,0)_{p,q}$ links
to $(q-1,p-1-i)_{q,p},$ which lies in the interval corresponding
to the vertex $v_{\frac{pq}{d}-1}$.
 Therefore, all the vertices
$v_{\alpha ((i,0)_{p,q})}, i\in [0,d-1],$ are from a common $\sim$
equivalence class.
 Since the elements in the interval $I _i$ are
evidently all $\approx$ equivalent to  $(i,0)_{p,q},$ we find that
$ \cup _{i=0}^{d-1}I_i$ belongs to one $\approx$ equivalence
class, say $A$. But the fact that $d\leq q$ implies that $I _i\cap
I _j=\emptyset$ as long as $i\not= j.$ Consequently, we deduce
from $\mid A\mid =d^2$ that $\mid A\mid = \mid \cup
_{i=0}^{d-1}I_i\mid $ and henceforth $A=\cup _{i=0}^{d-1}I_i.$
From $p>d,$ we see that $I_{p-1}\cap A=\emptyset.$ In particular,
this gives $(p-1,0)_{p,q}\notin A,$ and hence $(0,1)_{p,q}\notin
A,$ in virtue of Eq. (\ref{final}). This is impossible as  we
surely have $(0,1)_{p,q}\in I_0\subseteq A.$ This is the end of
the proof.
\end{proof}

It is immediate from Theorem \ref{theo} that Conjecture \ref{conj}
is true for $n=1.$ Hence,  we know that there is a unique OTIS
layout for $B(d,1)=K_d^+.$

\bibliography{pom}
\end{document}